 \documentclass[12pt,preprint]{aastex}

\newcommand{\mathsym}[1]{{}}
\newcommand{\unicode}[1]{{}}
 \usepackage{graphics}
 \usepackage{graphicx}
 \usepackage{epsfig}

\usepackage{amssymb}
\usepackage{amsmath}
\usepackage{aas_macros}

\newcommand{\lgbya}{\log\frac{b}{a}}

\shorttitle{Analysis of Linear Stability\dots Power-law Profile}
\shortauthors{Badam Singh Kushvah et al.}

\begin{document}


\title{Existence of Equilibrium Points and their Linear Stability in the Generalized Photogravitational Chermnykh-Like Problem with Power-law Profile}


\author{Badam Singh Kushvah,   Ram Kishor}
\affil{Department of   Applied  Mathematics,
Indian School of Mines, Dhanbad - 826004, Jharkhand,India}
%
\email{bskush@gmail.com;  kishor.ram888@gmail.com} 
 \and
\author{Uday Dolas}
 \affil{Department Of Mathematics, Chandra Shekhar Azad Govt. P.G. Nodal College,Sehore- 466001,  M.P., India}
\email{udolas@gmail.com}



\begin{abstract}
We consider the modified restricted three body problem with power-law density profile of disk, which rotates around the center of mass of the system with perturbed mean motion. Using analytical and numerical methods we have found equilibrium points and examined their linear stability. We have also found the zero velocity surfaces for the present model. In addition to five equilibrium points there is a new equilibrium point on the line joining the two primaries. It is found that $L_2$ and $L_3$ are stable for some values of inner and outer radius of the disk while collinear points are unstable, but $L_4$ is conditionally stable for mass ratio less than that of Routh's critical value. Lastly we have obtained the effects of radiation pressure, oblateness and mass of the disk.\end{abstract}





\section{Introduction}
The problem, after imposing a restriction as one body of the three body problem is of an infinitesimal(negligible) mass and remaining  other two are of finite masses, is known as restricted three body problem(RTBP). The governing force of motion of the RTBP is mainly the gravitational forces exerted by the finite masses also known as primaries. In the RTBP if we take bigger primary as a source of radiation then problem called as photogravitational RTBP which is generalized by taking smaller primary as an oblate spheroid. 
 
The Chermnykh-like problem which was first time studied by \cite{Chermnykh1987VeLen.......73C}, deals the motion of an infinitesimal mass in the orbital plane of a disk which rotates around the center of mass of the primaries with constant angular velocity $n$. \cite{K.Gozdzieski1998CeMDA..70...41G}, examined the problem in the sense of nonlinear stability of equilibrium points and also obtained the range of parameter for the same.

 The Chermnykh-Like problem has a number of applications in different areas as celestial mechanics (\cite{Gozdziewski1999CeMDA..75..251G}), chemistry (\cite{Strand1979JChPh..70.3812S}) etc. Also, the importance of the problem have seen in the extra solar planetary system (\cite{Rivera2000ApJ...530..454R}, \cite{Jiang2001AA...367..943J} etc.). 

 Further the effect of disc is very helpful in the study of resonance capture of Kuiper Belt Objects (KBOs) as given in \cite{Jiang2004MNRAS.355L..29J}. \cite{Papadakis2005Ap&SS.299..129P}, by taking assumptions as constant mass parameter and variable angular velocity parameter, analyzed the equilibrium point and zero velocity curve; \cite{Papadakis2005Ap&SS.299...67P} also studied the problem numerically. \cite{Jiang2006Ap&SS.305..341J} examined the Chermnykh problem with {$\mu$} = 0.5 and shown a deviation in the result of classical RTBP; also they have found the new equilibrium points in spite of Lagrangian points.  \cite{YehLCJiang2006} have found the condition of existence of new equilibrium points analytically and numerically.

 \cite{IshwarKushvah2006math......2467I} examined the linear stability of triangular points with P-R drag. Again \cite{Kushvah2008Ap&SS.318...41K} examined the stability of collinear points and found  unstable points.

Motivating by the importance and applications of the Chermnykh-like problem, we modeled a generalized photogravitational Chermnykh-like problem (section-\ref{sec:model}) in which we consider the angular velocity parameter $n>1$ which depends on the gravitational force of the disk assumed, radiation force of radiating (bigger) primary and oblateness factor of smaller primary which is an oblate spheroid.Further we determine the equilibrium points (section-\ref{sec:equipts}) and zero velocity surface (section-\ref{sec:zvs}) and then find the stability of the equilibrium points (section-\ref{sec:stb}) finally conclude the results (section-\ref{sec:conclude}).

\section{Mathematical Formulation of Model}
\label{sec:model}
Let us consider the motion of an infinitesimal mass governed by the gravitational force from the radiating body of mass $m_{1}$, oblate spheroid of mass $m_{2}$ ($m_1>m_2$) and a disk which considered, around the central binary system, of thickness $h\approx{10^{-4}}$ with power-law density profile $\rho(r)=\frac{c}{r^{p}}$, where $p$ is natural number(here we take $p=3$) and $c$ is a constant determined by the help of total mass of the disk.

We have chosen the unit of mass such that $G(m_1+m_2)={1}$; the primaries are separated by unit distance so that unit of time obtained by the choice taken; $\mu=\frac{m_2}{m_1+m_2}$  be the mass parameter then we have the masses $Gm_1=\mu$ and $Gm_2=1-\mu$.

Let us suppose $Oxyz$ be the rotating coordinate system having origin $O$ at the center of mass of the primaries which is fixed with respect to the inertial system and angular velocity $\omega$ along $z$-axis. $P(x,y,0)$, $A(x_1,0,0)=(-\mu,0,0)$ and $B(x_2,0,0)=(1-\mu,0,0)$ be the positions of infinitesimal body, radiating and oblate primaries respectively, relative to the rotating system.

 So, the equations of motion of the infinitesimal mass in $xy$-plane with respect to assumptions taken above be written as[as in \cite{Kushvah2011Ap&SS.332...99K}]

\begin{eqnarray}
\ddot{x}-2n\dot{y}&=&\Omega_{x}, \label{eq:ux}\\
\ddot{y}+2n\dot{x}&=&\Omega_{y}, \label{eq:vx}  \end{eqnarray}
where\begin{eqnarray*}
 &&\Omega_x=n^{2}x-\frac{(1-\mu)q_1(x+\mu)}{r^{3}_1}-\frac{\mu(x+\mu-1)}{r^{3}_2}-\frac{3}{2}\frac{\mu{A_2}(x+\mu-1)}{r^{5}_2}-V_{x}, \nonumber\\
&&\Omega_y=n^{2}y-\frac{(1-\mu)q_1{y}}{r^{3}_1}-\frac{\mu{y}}{r^{3}_2}-\frac{3}{2}\frac{\mu{A_2}{y}}{r^{5}_2}-V_{y}, \nonumber\ \end{eqnarray*}
and\begin{eqnarray}
&&\Omega=\frac{n^{2}(x^{2}+y^{2}}{2}+\frac{(1-\mu)q_1}{r_1}+\frac{\mu}{r_2}+\frac{\mu{A_2}}{{2}r^{3}_2}-V, \end{eqnarray}

with $r_1=\sqrt{(x+\mu)^{2}+y^{2}}$, $r_2=\sqrt{(x+\mu-1)^{2}+y^{2}}$ and  $r=\sqrt{x^{2}+y^{2}}$.

In above expression $q_1=(1-\beta)=(1-\frac{F_p}{F_g})$, the mass reduction factor of radiating body, $F_p$ is the radiation pressure force and $F_g$ is the gravitational force of same primary.
$A_2=\frac{R^{2}_e-R^{2}_p}{5R^{2}}$ is the  oblateness coefficient of second primary[as in \cite{McCuskey1963QB351.M3}], $R_e$ and $R_p$ are the equatorial and polar radii of the same body respectively and $R$ is the distance between primaries. Then the  mean motion of the primaries is given as 
$n=\sqrt{q_1+\frac{3}{2}A_2-2f_b(r)}$ which is greater than $1$.
The potential $(V)$ and gravitational force $f_b(r)$, of the disk are given as:
\begin{eqnarray}
&&V=-2 c h \pi \frac{(b-a)}{ab}\frac{1}{r}+\frac{7}{8} c h \pi \frac{\lgbya}{r^{2}} , \\
&&f_b(r)=-2 c h \pi \frac{(b-a)}{ab}\frac{1}{r^{2}}-\frac{3}{8} c h \pi \frac{\lgbya}{r^{3}}, \end{eqnarray}
where, a and b are inner and outer radii of the disk respectively. We assume that the gravitational force $f_b(r)$ is radially symmetric, so we have $\frac{x}{r}f_b(r)$ and $\frac{y}{r}f_b(r)$ as $x$ and $y$ components of the force $f_b(r)$ respectively.
Now, with the help of equations (\ref{eq:ux}) and (\ref{eq:vx}) we find the energy integral of the problem,
\begin{eqnarray}
&&C=-\dot{x}^{2}-\dot{y}^{2}+2\Omega, \label{eq:ji} \end{eqnarray}
where constant C is known as Jacobi constant.

\section{Equilibrium Points}
\label{sec:equipts}
The coordinates of equilibrium points of the problem are  obtained[as in \cite{Moulton1914QB351.M92}] by equating R.H.S.of the equations (\ref{eq:ux}) and (\ref{eq:vx}) both, to zero i.e. $\Omega_{x}=\Omega_{y}=0$ and solving them for $x$ and $y$. In other words \begin{eqnarray}
&& n^{2}x-\frac{(1-\mu)q_1(x+\mu)}{r^3_1}-\frac{\mu(x+\mu-1)}{r^3_2}-\frac{3}{2}\frac{\mu{A_2}(x+\mu-1)}{r^5_2}\nonumber\\&&-2 c h \pi \frac{(b-a)}{ab}\frac{x}{r^{3}}-\frac{3}{8} c h \pi {\lgbya}\frac{x}{r^{4}}=0 \label{eq:wx} \\
&&n^{2}y-\frac{(1-\mu)q_{1}{y}}{r^3_1}-\frac{\mu{y}}{r^3_2}-\frac{3}{2}\frac{\mu{A_2}y}{r^5_2}-2 c h \pi \frac{(b-a)}{ab}\frac{y}{r^{3}}-\nonumber\\&&\frac{3}{8} c h \pi {\lgbya}\frac{y}{r^{4}}=0 \label{eq:wy} \end{eqnarray}
Solving equations (\ref{eq:wx}) and (\ref{eq:wy}) for $x$ and $y$ we get the equilibrium points, separately given in the following subsections.
\subsection{Collinear Equilibrium Points}
For the collinear points we have $y=0$, so that $r_{1}=|x+\mu|$, $r_{2}=|x+\mu-1|$ and $r=|x|$.
Now suppose \begin{eqnarray}
&&f(x,y)=n^{2}x-\frac{(1-\mu)q_1(x+\mu)}{r^3_1}-\frac{\mu(x+\mu-1)}{r^3_2}-\frac{3}{2}\frac{\mu{A_2}(x+\mu-1)}{r^5_2} \nonumber\\&&-2 c h \pi \frac{(b-a)}{ab}\frac{x}{r^{3}}-\frac{3}{8} c h \pi {\lgbya}\frac{x}{r^{4}}=0, \label{eq:fx}\\
&&g(x,y)=n^{2}y-\frac{(1-\mu)q_{1}{y}}{r^3_1}-\frac{\mu{y}}{r^3_2}-\frac{3}{2}\frac{\mu{A_2}y}{r^5_2}-2 c h \pi \frac{(b-a)}{ab}\frac{y}{r^{3}}  \nonumber\\&&-\frac{3}{8} c h \pi {\lgbya}\frac{y}{r^{4}}=0. \label{eq:gy} \end{eqnarray} 
Consequently by substituting $y=0$ in equation (\ref{eq:fx}), we have $f(x,0)=0=K(x)$ (say) i.e. 
 \begin{eqnarray}
 &&K(x)=n^{2}x-\frac{(1-\mu)q_1(x+\mu)}{|x+\mu|^3}-\frac{\mu(x+\mu-1)}{|x+\mu-1|^3}-\frac{3}{2}\frac{\mu{A_2}(x+\mu-1)}{|x+\mu-1|^5}\nonumber\\&&-2 {c h \pi \frac{(b-a)}{ab}\frac{x}{|x|^3}}-\frac{3}{8} c h \pi \lgbya\frac{x}{|x|^4}=0. \label{eq:kx} \end{eqnarray}

Now, for the sake of simplicity we first divide the plane of motion $Oxy$ into three parts relative to the primaries as  $1-\mu\leq x$ , $-\mu<x<1-\mu$  and $x\leq-\mu$, we further divide second part into two sub parts as $0\leq x<1-\mu$ , $-\mu<x<0$, and for each considerable interval we have the function $K(x)$, which we shall use in equation (\ref{eq:kx}) for further analysis, as follows:\begin{equation}
K(x)=\begin{cases} n^{2}x-\frac{(1-\mu)q_1}{(x+\mu)^2}-\frac{\mu}{(x+\mu-1)^2}-\frac{3}{2}\frac{\mu{A_2}}{(x+\mu-1)^4}\\-2 c h \pi \frac{(b-a)}{ab}\frac{1}{x^2}-\frac{3}{8} c h \pi {\lgbya}\frac{1}{x^3}  & \quad \text{If}\  1-\mu<x,\\
n^{2}x-\frac{(1-\mu)q_1}{(x+\mu)^2}+\frac{\mu}{(x+\mu-1)^2}+\frac{3}{2}\frac{\mu{A_2}}{(x+\mu-1)^4}\\-2 c h \pi \frac{(b-a)}{ab}\frac{1}{x^2}-\frac{3}{8} c h \pi {\lgbya}\frac{1}{x^3}   &\quad \text{If}\ 0\leq x<1-\mu,\\
n^{2}x-\frac{(1-\mu)q_1}{(x+\mu)^2}+\frac{\mu}{(x+\mu-1)^2}+\frac{3}{2}\frac{\mu{A_2}}{(x+\mu-1)^4}\\+2 c h \pi \frac{(b-a)}{ab}\frac{1}{x^2}-\frac{3}{8} c h \pi {\lgbya}\frac{1}{x^3}   &\quad \text{If}\ -\mu<x<0,\\
n^{2}x+\frac{(1-\mu)q_1}{(x+\mu)^2}+\frac{\mu}{(x+\mu-1)^2}+\frac{3}{2}\frac{\mu{A_2}}{(x+\mu-1)^4}\\+2 c h \pi \frac{(b-a)}{ab}\frac{1}{x^2}-\frac{3}{8} c h \pi {\lgbya}\frac{1}{x^3}   &\quad \text{If}\ x<-\mu,       \end{cases} \end{equation}

Case(1)  when $1-\mu<x$: Let the distance between $\mu$ and equilibrium point on the $x$-axis in the interval $\left[ 1-\mu,+\infty\right)$ be $\rho$ then $x-x_2=(x+\mu-1)=\rho>0$  and so in this case $x-x_1=(x+\mu)=(1+\rho)>0$ and $x=(1+\rho-\mu)>0$. Substituting these values in equation (\ref{eq:kx}) and simplifying, we have 
\begin{eqnarray}&&8 n^2 \rho ^{10}+B_1 \rho ^9+B_2 \rho ^8+ \left(B_3-16 c h \pi \frac{(b-a)}{ab}-8 \mu \right) \rho ^7+\{B_4-48 c h \pi \frac{(b-a)}{ab}+\nonumber\\&&8\left(2 c h \pi \frac{(b-a)}{ab}-5\right) \mu+24 \mu ^2\}\rho ^6+\{B_5-48 c h \pi \frac{(b-a)}{ab} +16\left(2 c h \pi \frac{(b-a)}{ab}-5\right) \mu+\nonumber\\&&96 \mu ^2-24 \mu ^3\}\rho ^5+\{B_6-16 c h \pi \frac{(b-a)}{ab}+16\left(c h \pi \frac{(b-a)}{ab}-5\right) \mu+144 \mu ^2-72 \mu ^3+8 \mu ^4\} \rho ^4\nonumber\\&&+B_7 \rho ^3+B_8 \rho ^2+B_9 \rho  +B_{10}=0. \label{eq:rx}  \end{eqnarray}
where, $B_1= 48 n^2-32 n^2 \mu$, $B_2= 120 n^2-160 n^2 \mu +48 n^2 \mu ^2$,

$B_3= 160 n^2-8q_1+ 8\left(-40 n^2+ q_1\right) \mu +192 n^2 \mu ^2-32 n^2 \mu ^3$,

$B_4=120 n^2-3 c h \pi {\lgbya}-24 q_1+16\left(3 q_1-20 n^2\right)\mu+24\left(-q_1+12 n^2\right)\mu ^2-96\mu ^3+8 n^2\mu ^4$,

\begin{eqnarray}&&B_5=48 n^2-6 c h \pi {\lgbya}-24 q_1-4\left(40 n^2+3 A_2-18 q_1\right)\mu+24\left(8 n^2-3 q_1\right)\mu ^2+\nonumber\\&&24\left(-4 n^2+q_1\right)\mu ^3+16 n^2 \mu ^4,\nonumber \end{eqnarray}
\begin{eqnarray}&&B_6=8 n^2-3 c h \pi {\lgbya}-8 q_1-4\left(8 n^2+15 A_2-8 q_1\right)\mu+12\left(4 n^2+3 A_2-4 q_1\right)\mu ^2+\nonumber\\&&32\left(-n^2+q_1\right)\mu ^3+8\left(n^2-q_1\right)\mu ^4,\nonumber \end{eqnarray}

$B_7=-40\left(1+3 A_2\right)\mu+48\left(2+3 A_2\right)\mu ^2-36\left(2+A_2\right)\mu ^3+16 \mu ^4$,

$B_8=-8\left(1+15 A_2\right)\mu+24\left(1+9 A_2\right)\mu ^2-12\left(2+9 A_2\right)\mu ^3+4\left(2+3 A_2\right)\mu ^4$,

$B_9=-60 A_2 \mu+144 A_2 \mu ^2-108 A_2 \mu ^3+24 A_2 \mu ^4$ and $B_{10}=-12 A_2 \mu+36 A_2 \mu ^2-36 A_2 \mu ^3+12 A_2 \mu ^4$.

Suppose all the quantities are constants except $\mu$ on which the roots of the above equation depend. So, here we assume that L.H.S. of the equation (\ref{eq:rx}) is the function of $\rho$ and $\mu$.
For $\mu=0$, we have 
\begin{eqnarray}&&\rho ^4\left(8 n^2 \rho ^6+48 n^2 \rho ^5+120 n^2 \rho ^4+C_1 \rho ^3+C_2 \rho ^2+C_3 \rho+C_4\right)=0, \label{eq:rt} \end{eqnarray}
where, $C_1=160 n^2-16 c h \pi \frac{(b-a)}{ab}-8 q_1$, $C_2=120 n^2-48 c h \pi \frac{(b-a)}{ab}-3 c h \pi {\lgbya}-24 q_1$,

$C_3=48 n^2-48 c h \pi \frac{(b-a)}{ab}-6 c h \pi {\lgbya}-24 q_1$ and $C_4=8 n^2-16 c h \pi \frac{(b-a)}{ab}-3 c h \pi {\lgbya}-8 q_1$.

One can see that the equation (\ref{eq:rt}) has four roots equal to zero and others come from remaining factor. So applying the theory of solution of the algebraic equation assuming as $\mu$ is very small, the four roots of the equation (\ref{eq:rt}) are expressible as a power series in $\mu^{1/4}$ and vanish with $\mu$. Two of the four roots are real and others are complex with the real value of $\mu^{1/4}$. There fore the series is given as:
\begin{eqnarray}&&\rho=\alpha_1 \mu ^{\frac{1}{4}}+\alpha_2 \mu ^{\frac{2}{4}}+\alpha_3 \mu ^{\frac{3}{4}}+\alpha_4 \mu ^{\frac{4}{4}}+\dots, \label{eq:ps}\end{eqnarray}

where $\alpha_1,\alpha_2,\alpha_3,\alpha_4,\dotsc$ are constant to be determined.
 Putting this $\rho$ into the equation (\ref{eq:rx}) and equating, the coefficient of corresponding power of  $\mu^{1/4}$, to the zero we get

$\alpha_1=\pm \frac{\alpha ^{\frac{1}{4}}}{d_1}$,  $\alpha_2={\frac{d_2}{4 d_1}}{\alpha ^{\frac{2}{4}}}$, $\alpha_3=\pm{\frac{\{d_1^2+3 A_2 (d_3+d_4)\}}{6 A_2 d_1^2}}{\alpha ^{\frac{3}{4}}}$ and 

$\alpha_4=-{\frac{\{2 d_1^2 d_2+3 a b A_2 (d_5+d_6)\}}{12 A_2 d_1^3}}{\alpha ^{\frac{4}{4}}}$, where 

$\alpha= \frac{12 a b A_2}{d_1}$, $d_1= 8 a b n^2-16 c h \pi(b-a)-3 a b c h \pi {\lgbya}-8 a b q_1$,

$d_2= 8 a b n^2+32 c h \pi (b-a)+9 a b c h \pi {\lgbya}+16 a b q_1$,
\begin{eqnarray}&&d_3= [448 a^2 b^2 n^4-6656 a^2 n^2 c h \pi (b-a)+1024 c^2 h^2 \pi ^2 (b-a)^2+135 a^2 b^2 c^2 h^2 \pi ^2 ({\lgbya})^2+\nonumber\\&&144 a b c h\{15 a b n^2+4 c h \pi (b-a)\}{\lgbya}],\nonumber \end{eqnarray}

$d_4=32 a b\{104 a b n^2+32 c h \pi (b-a)+9 a b q_1 c h \pi {\lgbya}\}+256 a^2 b^2 q_1^2$,
\begin{eqnarray}&&d_5= [512 a^2 b^2 n^6-16384 a^2 n^4 c h \pi (b-a) +20480 n^2 c^2 h^2 \pi ^2 (b-a)^2+5952 a^2 b^2 n^4 c h \pi {\lgbya}+\nonumber\\&&2232 a^2 b^2 n^2 c^2 h^2 \pi ^2 ({\lgbya})^2-27 a^2 b^2 c^3 h^3 \pi ^3 ({\lgbya})^3+12288 a b n^2 c^2 h^2 \pi ^2 (b-a){\lgbya}],\nonumber   \end{eqnarray}

$d_6=2048 a b n^2 \{4 a b n^2+10 b c h \pi (b-a)+3 a b c h \pi {\lgbya}\}q_1+5120 a^2 b^2 n^2 q_1^2$.

So after calculating these constants and with help of equation (\ref{eq:ps}) we get $\rho$. Therefore in the case (1) we have \begin{equation}
\begin{cases}
 r_1=|x+\mu|=1+\rho=1+\alpha_1 \mu ^{\frac{1}{4}}+\alpha_2 \mu ^{\frac{2}{4}}+\alpha_3 \mu ^{\frac{3}{4}}+\alpha_4 \mu ^{\frac{4}{4}}+\dots,  \\r_2=|x+\mu-1|=\rho=\alpha_1 \mu ^{\frac{1}{4}}+\alpha_2 \mu ^{\frac{2}{4}}+\alpha_3 \mu ^{\frac{3}{4}}+\alpha_4 \mu ^{\frac{4}{4}}+\dots, \\r=|x|=1- \mu+\rho=1-\mu+\alpha_1 \mu ^{\frac{1}{4}}+\alpha_2 \mu ^{\frac{2}{4}}+\alpha_3 \mu ^{\frac{3}{4}}+\alpha_4 \mu ^{\frac{4}{4}}+\dots.
\end{cases} \label{eq:zx1} \end{equation}

Case(2)  when $-\mu<x<1-\mu$:   sub case(i)  when $0\leq x<1-\mu$:  let us suppose that the distance between $\mu$ and equilibrium point on the $x$-axis in the interval $\left[0,1-\mu\right)$ be $-\rho$ then $x-x_2=(x+\mu-1)=-\rho<0$  and so in this case $x-x_1=(x+\mu)=(1-\rho)>0$ and $x=(1-\rho-\mu)>0$. Putting these values in equation (~\ref{eq:kx}) for same interval and simplifying, we have changed form  of equation (\ref{eq:rx}) as follows: 

\begin{eqnarray}&&8 n^2 \rho ^{10}-B_1 \rho ^9+B_2 \rho ^8+ \left(-B_3+16 c h \pi \frac{(b-a)}{ab}-8 \mu \right) \rho ^7+\{B_4-48 c h \pi \frac{(b-a)}{ab}+\nonumber\\&&8\left(16 c h \pi  \frac{(b-a)}{ab}+5\right) \mu-24 \mu ^2\} \rho ^6+\{-B_5+48 c h \pi \frac{(b-a)}{ab} -16\left(2 c h \pi  \frac{(b-a)}{ab}+5\right) \mu+\nonumber\\&&96 \mu ^2-24 \mu ^3\}\rho ^5+\{B_6-16 c h \pi \frac{(b-a)}{ab}+16\left(16 c h \pi \frac{(b-a)}{ab}-5\right) \mu-144 \mu ^2+72 \mu ^3-8 \mu ^4\} \rho ^4\nonumber\\&&-B_7 \rho ^3+B_8 \rho ^2-B_9 \rho  +B_{10}=0, \label{eq:rx2i} \end{eqnarray}

also changed form  of equation (\ref{eq:rt}) is 
\begin{eqnarray}&&\rho ^4\left(8 n^2 \rho ^6-48 n^2 \rho ^5+120 n^2 \rho ^4-C_1 \rho ^3+C_2 \rho ^2-C_3 \rho+C_4\right)=0,  \end{eqnarray} where $B_1, B_2,\dots, B_9, B_{10}, C_1, C_2, C_3$ and $C_4$ are given as in case (1).
 Similar analysis as in case (1) provides  
\begin{eqnarray}\begin{cases}
 r_1=|x+\mu|=1-\rho=1-\alpha_1 \mu ^{\frac{1}{4}}-\alpha_2 \mu^{\frac{2}{4}}-\alpha_3 \mu^{\frac{3}{4}}-\alpha_4 \mu ^{\frac{4}{4}}-\dots, \\r_2=|x+\mu-1|=-\rho=-\alpha_1 \mu ^{\frac{1}{4}}-\alpha_2 \mu ^{\frac{2}{4}}-\alpha_3 \mu ^{\frac{3}{4}}-\alpha_4 \mu^{\frac{4}{4}}-\dots, \\r=|x|=1- \mu-\rho=1-\mu-\alpha_1 \mu^{\frac{1}{4}}-\alpha_2 \mu ^{\frac{2}{4}}-\alpha_3 \mu ^{\frac{3}{4}}-\alpha_4 \mu ^{\frac{4}{4}}-\dots.
\end{cases} \label{eq:zx2i}\end{eqnarray} where constant $\alpha_1,\alpha_2,\alpha_3,\alpha_4,\dots $ are given as:

$\alpha_1=\pm \frac{\left(-\alpha\right) ^{\frac{1}{4}}}{d_1}$,  $\alpha_2={\frac{d_2}{4 d_1}}{\left(-\alpha\right) ^{\frac{2}{4}}}$, $\alpha_3=\pm{\frac{\{d_1^2+3 A_2 (d_3+d_4)\}}{6 A_2 d_1^2}}{\left(-\alpha\right) ^{\frac{3}{4}}}$ and 

$\alpha_4=-{\frac{\{2 d_1^2 d_2+3 a b A_2 (d_5+d_6)\}}{12 A_2 d_1^3}}{\left(-\alpha\right) ^{\frac{4}{4}}}$ and $\alpha, d_1, d_2, d_3, d_4, d_5$ and $d_6$ are as in case (1).

Sub case(ii)  when $-\mu<x<0$: let the distance between $\mu$ and equilibrium point on the $x$-axis in the interval $\left(-\mu, 0\right)$ be $-\rho$ then $x-x_2=(x+\mu-1)=-\rho<0$, $x-x_1=(x+\mu)=(1-\rho)>0$ and $x=(1-\rho-\mu)<0$ (The main difference in sub case (i) and(ii) is that in (i) $x>0$ and in (ii) $x<0$). Putting these values in equation (\ref{eq:kx}) for same interval and simplifying, we have changed form  of equation (\ref{eq:rx}) of the case (1) as follows:
 \begin{eqnarray}&&8 n^2 \rho ^{10}-B_1 \rho ^9+B_2 \rho ^8+ \left(-B_3-16 c h \pi \frac{(b-a)}{ab}-8 \mu \right) \rho ^7+\{B_4+48 c h \pi \frac{(b-a)}{ab}-\nonumber\\&&8\left(2 c h \pi \frac{(b-a)}{ab} -5\right) \mu-24 \mu ^2\}\rho ^6+\{-B_5-48 c h \pi \frac{(b-a)}{ab} +16\left(2 c h \pi \frac{(b-a)}{ab}-5\right) \mu+\nonumber\\&&96 \mu ^2-24 \mu ^3\}\rho ^5+\{B_6+16 c h \pi \frac{(b-a)}{ab}-16\left(c h \pi \frac{(b-a)}{ab}-5\right) \mu-144 \mu ^2+72 \mu ^3-8 \mu ^4\} \rho ^4-\nonumber\\&&B_7 \rho ^3+B_8 \rho ^2-B_9 \rho  +B_{10}=0, \label{eq:rx2ii} \end{eqnarray}  also, changed form  of equation (\ref{eq:rt}) of the case (1) is \begin{eqnarray}
&&\rho ^4\biggl[8 n^2 \rho ^6-48 n^2 \rho ^5+120 n^2 \rho ^4+\left(C_1+16 q_1\right) \rho ^3+\left(C_2+96 c h \pi \frac{(b-a)}{ab}\right) \rho ^2\nonumber\\&&+\left(C_3-12 c h \pi {\lgbya}+48 q_1\right) \rho+\left(C_4+32 c h \pi \frac{(b-a)}{ab}\right)\biggr]=0, \end{eqnarray} where $B_1, B_2,\dots, B_9, B_{10}, C_1, C_2, C_3$ and $C_4$ are given as in case (1).
\begin{eqnarray} \begin{cases}
 r_1=|x+\mu|=1-\rho=1-\alpha_1 \mu ^{\frac{1}{4}}-\alpha_2 \mu ^{\frac{2}{4}}-\alpha_3 \mu ^{\frac{3}{4}}-\alpha_4 \mu ^{\frac{4}{4}}-\dots, \\r_2=|x+\mu-1|=-\rho=-\alpha_1 \mu ^{\frac{1}{4}}-\alpha_2 \mu ^{\frac{2}{4}}-\alpha_3 \mu ^{\frac{3}{4}}-\alpha_4 \mu ^{\frac{4}{4}}-\dots, \\r=|x|=-(1- \mu-\rho)=-1+\mu+\alpha_1 \mu ^{\frac{1}{4}}+\alpha_2 \mu ^{\frac{2}{4}}+\alpha_3 \mu ^{\frac{3}{4}}+\alpha_4 \mu ^{\frac{4}{4}}+\dots.\end{cases} \label{eq:zx2ii} \end{eqnarray} where constant $\alpha_1,\alpha_2,\alpha_3,\alpha_4,\dotsc$ and $\alpha$ given as in case (1) whereas

 $d_1= 8 a b n^2+16 c h \pi (b-a)-3 a b c h \pi {\lgbya}-8 a b q_1$,

$d_2= 8 a b n^2-32 a c h \pi (b-a) +9 a b c h \pi {\lgbya}+16 a b q_1$,
\begin{eqnarray}&&d_3= [448 a^2 b^2 n^4-6656 a b c h n^2 \pi (b-a)+1024 c^2 h^2 \pi ^2 (b-a)^2+135 a^2 b^2 c^2 h^2 \pi ^2 ({\lgbya})^2+\nonumber\\&&144 a b c h\{15 a b n^2-4 c h \pi (b-a)\}{\lgbya}],\nonumber  \end{eqnarray}

$d_4=32 a b \{104 a b n^2-32 c h \pi (b-a)+9 a b c h \pi {\lgbya}\}q_1+256 a^2 b^2 q_1^2$,
\begin{eqnarray}&&d_5= [512 a^2 b^2 n^6-16384 a b n^4 c h \pi (b-a)+20480 n^2 c^2 h^2 \pi ^2(b-a)^2+5952 a^2 b^2 n^4 c h \pi {\lgbya}+\nonumber\\&&2232 a^2 b^2 n^2 c^2 h^2 \pi ^2 ({\lgbya})^2+27 a^2 b^2 c^3 h^3 \pi ^3 ({\lgbya})^3-12288 a b n^2 c^2 h^2 \pi ^2 (b-a){\lgbya}],\nonumber   \end{eqnarray}
$d_6=2048 a b n^2 \{4 a b n^2-10 c h \pi (b-a)+3 a b c h \pi {\lgbya}\} q_1+5120 a^2 b^2 n^2 q_1^2$.

Case(3)  when $x<-\mu$: let the distance between $1-\mu$ and equilibrium point on the $x$-axis in the interval $\left(-\infty,-\mu\right)$  be denoted by $1-\rho$ then $x-x_1=\left(x+\mu\right)=\left(\rho-1\right)<0$ is negative and so $x-x_2=\left(x+\mu-1\right)=\left(\rho-2\right)<0$ and $x=\left(\rho-\mu-1\right)<0$. Putting these values in equation (\ref{eq:kx}) and simplifying we get 
\begin{eqnarray}
&&8 n^2 \rho ^{10}+B_1 \rho ^9+B_2 \rho ^8+B_3 \rho ^7+B_4 \rho ^6+B_5 \rho ^5+B_6 \rho ^4+B_7 \rho ^3+B_8 \rho ^2+B_9 \rho +B_{10}=0, \label{eq:rx3}  \end{eqnarray} where in this case,

$B_1=-112 n^2-32 n^2 \mu$, $B_2=696 n^2+416 n^2 \mu +48 n^2 \mu ^2$, 

$B_3=-2528 n^2+16 c h \pi \frac{(b-a)}{ab}+8 q_1-8\left(296 n^2+q_1-1\right)$,

\begin{eqnarray}&&B_4=5944 n^2-176 c h \pi \frac{(b-a)}{ab}-3 c h \pi {\lgbya}-88 q_1+8\{968 n^2-16 c h \pi \frac{(b-a)}{ab}\nonumber\\&&-8 q_1-9\}\mu+24\left(124 n^2+q_1-1\right)\mu ^2+352 n^2 \mu ^3+8 n^2 \mu ^4,\nonumber  \end{eqnarray}
\begin{eqnarray}&&B_5=-9456 n^2+816 c h \pi \frac{(b-a)}{ab}+30 c h \pi {\lgbya}+408 q_1+4\{-4008 n^2+40 c h \pi \frac{(b-a)}{ab}+\nonumber\\&&3 A_2-42 q_1+68\}\mu-8\left(1080 n^2+27 q_1-24\right)\mu ^2-24\left(68 n^2+q_1-1\right)\mu ^3-80 n^2 \mu ^4, \nonumber \end{eqnarray}
\begin{eqnarray}&&B_6=10312 n^2-2064 c h \pi \frac{(b-a)}{ab}-123 c h \pi {\lgbya}-1032 q_1+4\{5448 n^2-164 c h \pi \frac{(b-a)}{ab}\nonumber\\&&-15 A_2+12 q_1-140\}\mu+12\left(1284 n^2-3 A_2-64 q_1-52\right)\mu ^2+8\left(516 n^2+26 q_1-21\right)\mu ^3\nonumber\\&&+8\left(41 n^2+q_1-1\right)\mu ^4, \nonumber \end{eqnarray}
\begin{eqnarray}&&B_7=-7616 n^2+3072 c h \pi \frac{(b-a)}{ab}+264 c h \pi {\lgbya}+1536 q_1+8\{2432 n^2+176 c h \pi \frac{(b-a)}{ab}\nonumber\\&&+15 A_2+72 q_1+85\}\mu+4\left(-4320 n^2+36 A_2+336 q_1+264\right)\mu ^2+4\left(-1536 n^2+9 A_2-\right.\nonumber\\&&\left.176 q_1+114\right)\mu ^3-16\left(44 n^2+4 q_1-3\right)\mu ^4, \nonumber \end{eqnarray}
\begin{eqnarray}&&B_8=3648 n^2-2688 c h \pi \frac{(b-a)}{ab}-312 c h \pi {\lgbya}-1344 q_1+8\{1376 n^2-208 c h \pi \frac{(b-a)}{ab}\nonumber\\&&-15 A_2-144 q_1-61\}\mu+8\left(1488 n^2-27 A_2+144 q_1-123\right)\mu ^2+4\left(1344 n^2-27 A_2+288 q_1\right.\nonumber\\&&\left.-150\right)\mu ^3+4\left(208 n^2-3 A_2+48 q_1-13\right)\mu ^4, \nonumber  \end{eqnarray}
\begin{eqnarray}&&B_9=-1024 n^2+1280 c h \pi \frac{(b-a)}{ab}+192 c h \pi {\lgbya}+640 q_1+4\{-896 n^2+256 c h \pi \frac{(b-a)}{ab}\nonumber\\&&+15 A_2+224 q_1+48\}\mu+8\left(-576 n^2+18 A_2-48 q_1+60\right)\mu ^2+4\left(-640 n^2+27 A_2-224 q_1\right.\nonumber\\&&\left.+96\right)\mu ^3+8\left(-64 n^2+3 A_2-32 q_1+12\right)\mu ^4, \nonumber \end{eqnarray}
\begin{eqnarray}&&B_{10}=128 n^2-256 c h \pi \frac{(b-a)}{ab}-48 c h \pi {\lgbya}-128 q_1+4\{128 n^2-64 c h \pi \frac{(b-a)}{ab}\nonumber\\&&-3 A_2-64 q_1-8\}\mu+12\left(64 n^2-3 A_2-8\right)\mu ^2+4\left(128 n^2-9 A_2+64 q_1-24\right)\mu ^3\nonumber\\&&+4\left(32 n^2-3 A_2+32 q_1-8\right)\mu ^4. \nonumber \end{eqnarray} 
 
Putting $\mu=0$ in equation (\ref{eq:rx3}), we have \begin{eqnarray}
 &&(-2+\rho )^4 (-1+\rho )^2 [8 a b n^2 \rho ^4-32 a b n^2 \rho ^3+48 a b n^2 \rho ^2\{-32 a b n^2 +16 c h \pi (b-a)\nonumber\\&&+8 a b q_1\}\rho+8 a b n^2-16 c h \pi (b-a) -3 a b c h \pi {\lgbya}-8 a b q_1]=0.\label{eq:rt3} \end{eqnarray}
Clearly out of ten roots of equation (\ref{eq:rt3}), four equal to $2$, two equal to $1$ and others come from remaining factor. Now here we want to apply the theory of solution of the algebraic equation similar as earlier cases assuming as $\mu$ very small and L.H.S. of equation (\ref{eq:rx3}) is function of $\rho$ and $\mu$ only. So the four roots of the equation (\ref{eq:rx3}) are expressible as a power series in $\mu^{1/4}$ with an extra term as $2$ and two roots as a power series in $\mu^{1/2}$ with an extra term as $1$ and these series not vanish with $\mu$ due to extra term. The power series in $\mu^{1/2}$ with an extra term as $1$ is independent to $\mu$ as its constant coefficient comes out zero, so this series is unimportant whenever the power series in $\mu^{1/4}$  with an extra term as $2$ depends on $\mu$. The two of the four roots are real and other are complex with the real value of $\mu^{1/4}$.The power series is given as:
\begin{eqnarray}&&\rho=2+\alpha_1 \mu ^{\frac{1}{4}}+\alpha_2 \mu ^{\frac{2}{4}}+\alpha_3 \mu ^{\frac{3}{4}}+\alpha_4 \mu ^{\frac{4}{4}}+\dots, \label{eq:ps3} \end{eqnarray} 
where $\alpha_1,\alpha_2,\alpha_3,\alpha_4,\dotsc$ are given as in case (1) which is obtained by putting $\rho$ from (\ref{eq:ps3}) into the equation (\ref{eq:rx3}) and equating, the coefficient of corresponding power of  $\mu^{1/4}$, to the zero. In this case $\alpha$ is same as in earlier cases and $d_3$, $d_5$ are as in sub case (ii), whereas

$d_1= 8 a b n^2+16 c h \pi (b-a)-3 a b c h \pi {\lgbya}+8 a b q_1$,

$d_2= -8 a b n^2+32 c h \pi (b-a)-9 a b c h \pi {\lgbya}+16 a b q_1$,

$d_4=-32 a b \{104 a b n^2-32 c h \pi (b-a)+9 a b c h \pi {\lgbya}\}q_1+256 a^2 b^2 q_1^2$,

$d_6=-2048 a b n^2 \{4 a b n^2-10 c h \pi (b-a)+3 a b c h \pi {\lgbya}\} q_1+5120 a^2 b^2 n^2 q_1^2$.

Again from equation (\ref{eq:ps3}) we get $\rho$ and there fore in this case we have \begin{equation}
\begin{cases}
 r_1=|x+\mu|=-(\rho-1)=-1-\alpha_1 \mu ^{\frac{1}{4}}-\alpha_2 \mu ^{\frac{2}{4}}-\alpha_3 \mu ^{\frac{3}{4}}-\alpha_4 \mu ^{\frac{4}{4}}-\dots,  \\r_2=|x+\mu-1|=-(\rho-2)=-\alpha_1 \mu ^{\frac{1}{4}}-\alpha_2 \mu ^{\frac{2}{4}}-\alpha_3 \mu ^{\frac{3}{4}}-\alpha_4 \mu ^{\frac{4}{4}}-\dots, \\r=|x|=-(\rho- \mu-1)=-1+\mu-\alpha_1 \mu ^{\frac{1}{4}}-\alpha_2 \mu ^{\frac{2}{4}}-\alpha_3 \mu ^{\frac{3}{4}}-\alpha_4 \mu ^{\frac{4}{4}}-\dots   .\label{eq:zx3}
\end{cases} \end{equation}

For numerical calculation we set the parametric values as follows:
Mass parameter$(\mu)=\frac{m_j}{M_s+m_j}=0.000953728$, where $M_s$(Sun mass)$=1.98892\times10^{30}kg$ and $m_j$(Jupiter mass)$=1.8987\times10^{27}kg$, radiation factor$(q_1)=0.75$, oblateness effect$(A_2)=0.0025$, disk's mass$(M_b)=0.4$ (taking disk's inner radius$(a)=1$, disk's outer radius$(b)=1.5$, control factor of density profile$(c)=1910.83$, disk's thickness$(h)=0.0001$ and $\pi=3.14$). Putting these values into the equation (\ref{eq:kx}) for each intervals and solving we have collinear equilibrium points with respect to the relevant intervals given below:

In interval $(1-\mu, +\infty)$ we have six real values of $x$ out of that five lie outside the interval and remaining one i.e. $L_1=1.05667$ belongs to the interval. Similarly in case of interval $(0, 1-\mu)$ we have four real values of $x$, one of them i.e.$L_2=0.813609$ belongs to the interval $(0, 1-\mu)$ and three other lie outside the interval. In the interval, $(-\mu, 0)$ we again have four real values of $x$, one of them i.e.$l_1=-0.00879106$(new point) belongs to the interval $(-\mu, 0)$ and three other lie outside the interval but in the interval $(-\infty, -\mu)$, we have only two real value of $x$ one of them i.e. $L_3=-0.823420$ belongs to the $(-\infty, -\mu)$ and other one is positive. 

We also find the collinear equilibrium points with the help of equations (\ref{eq:zx1}), (\ref{eq:zx2i}), (\ref{eq:zx2ii}) and (\ref{eq:zx3}), after evaluating $\alpha_i$,$i=1,2,3,4\dots$ and hence $\rho$ for each case, also for mean motion $(n)$ we set here disk's reference radius $r=0.99$, given as $L_1=1.04411$, $L_2=0.971566$, $l_1=-0.975463$ and $L_3=-0.977926$ which is very close to the points obtained above one excepting the new point $l_1$ which is also outside of the relevant interval.
  \begin{figure}
  \plottwo{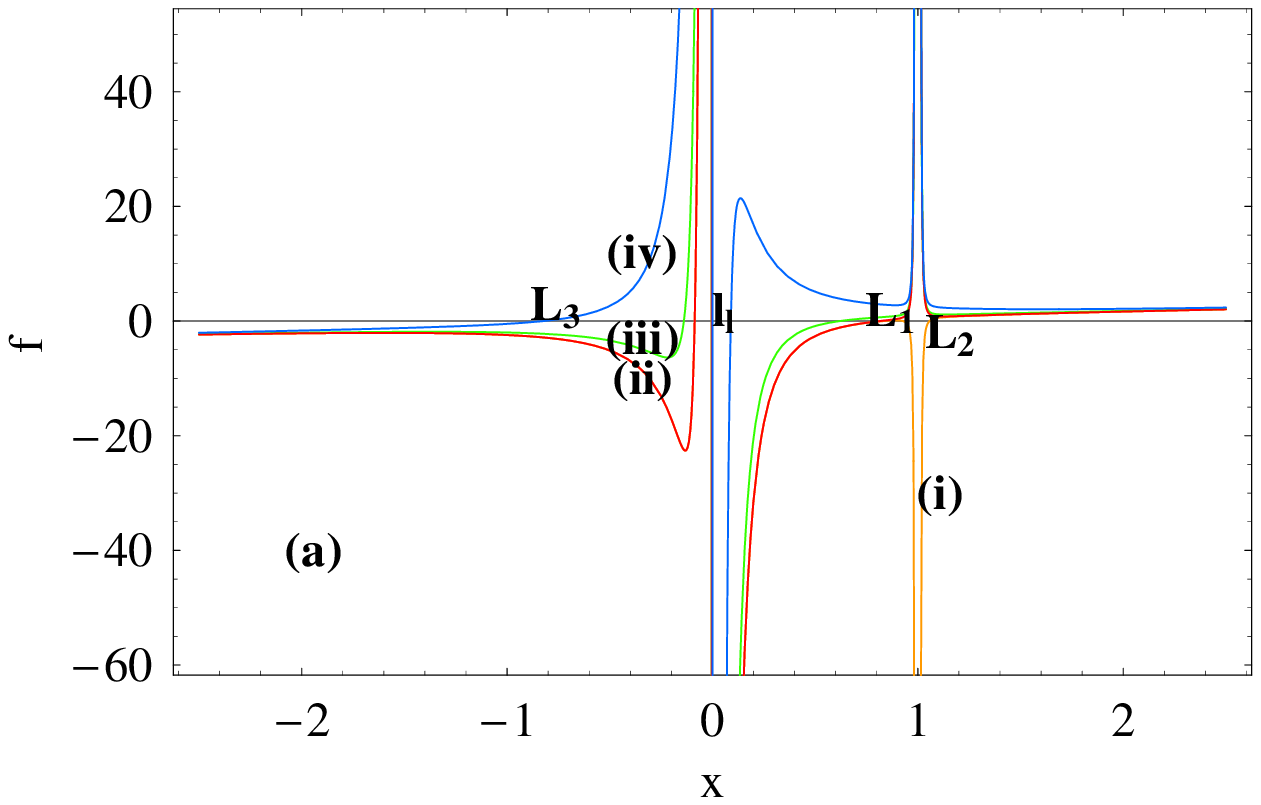}{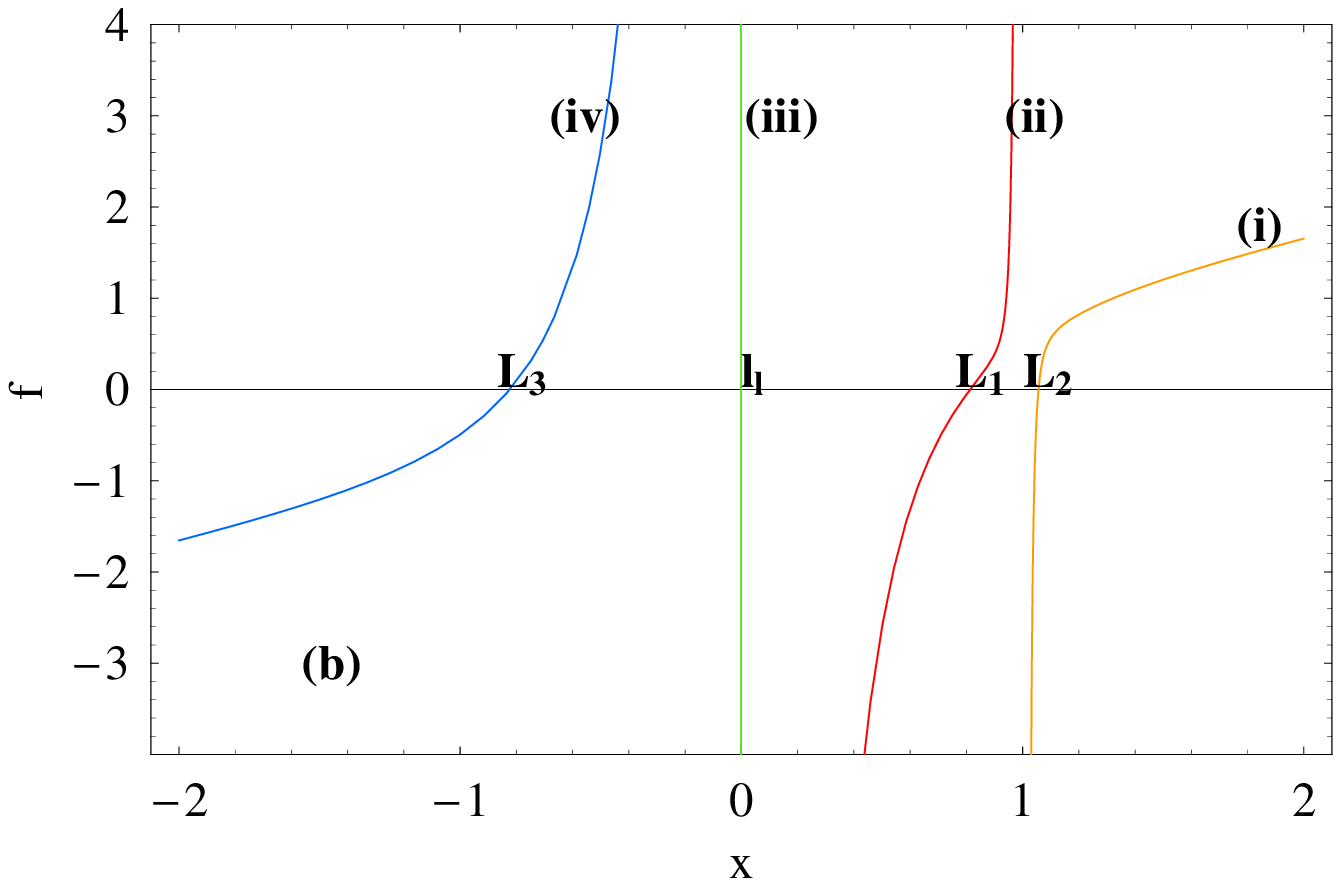}
  \caption{Collinear Points: $f(x, 0)$ Vs $x$- at  $q_1=0.75$, $A_2=0.0025$, $a=1$, $b=1.5$ (a) range and domain of $f(x, 0)$ both different (b)range and domain of $f(x, 0)$ both same.}\label{fa25q75b15lable}
 \end{figure}
   \begin{figure}
   \plotone{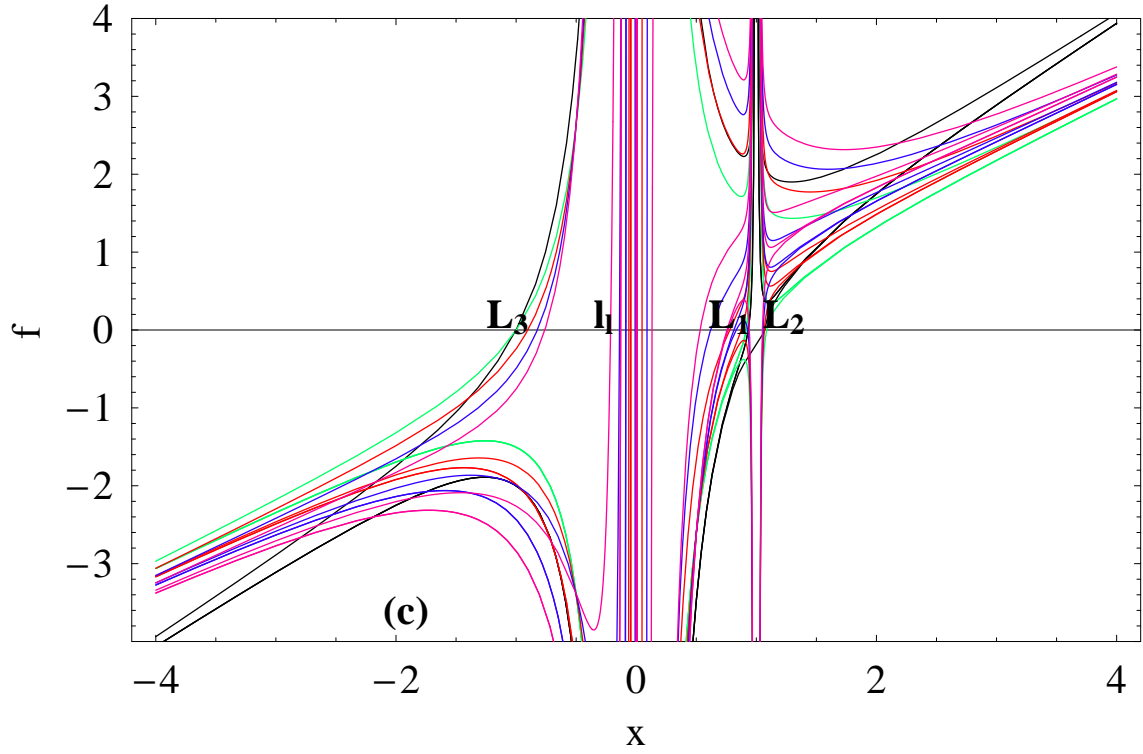}
  \caption{Comparison of nature of curve $f(x,0)$ drown at different values of parameters: $q_1, A_2$ and $b$.}\label{fq75a25lable}
 \end{figure}

In figure \ref{fa25q75b15lable}   $(a)$, $f(x, 0)=K(x)$ Vs $x$ contains four curves (i), (ii), (iii) and (iv) for the different intervals of $x$ or domains of $f(x, 0)$ but same range of $f(x, 0)$. These four curves intersect $x$-axis at one and only one point i.e. at $L_1$, $L_2$, $L_3$ and $l_1$ in their respective intervals (but out side the respective intervals we find some more intersections points of each curves which are other roots of the (\ref{eq:rx}), (\ref{eq:rx2i}), (\ref{eq:rx2ii}) and (\ref{eq:rx3}). These intersection points i.e. $L_1$, $L_2$, $L_3$ and $l_1$ are much clear in figure \ref{fa25q75b15lable}  $(b)$ which is drown at same domain and range, which is nothing but the respective intervals of each curves. In figure   \ref{fq75a25lable}, we plots same curves as in figure \ref{fa25q75b15lable} $(a)$ but at different values of disk's outer radius $b$ i.e at $b=1.0, 1.2, 1.5$, and $2.0$ in addition with classical case ($q_1=1$, $A_2=0$ and $a=b=1.0$). In other words figure   \ref{fq75a25lable}, shows the nature of curve $f(x, 0)$, at different widths of disk, which is similar in nature as in figure \ref{fa25q75b15lable}  $(a)$ and intersects $x$-axis at only one point $L_1$, $L_2$, $L_3$ and $l_1$  in their respective intervals while in classical case (even in our cases also i.e. if $q_1=0.75$, $A_2=0.0025$ and $a=b=1$) the curves $f(x, 0)$ intersect $x$-axis at three different points $L_1$, $L_2$, $L_3$ and for $x$ in $(-\mu, 0)$, $f(x, 0)$ does not intersect $x$-axis. 

Thus, we conclude that in the presence of disk, there are always four different intersection points $L_1$, $L_2$, $L_3$ and $l_1$ of curve $f(x, 0)$ with $x$-axis, each lies in the four different intervals described above.

\subsection{Triangular Equilibrium Point}
In this case we have $y\neq0$, for the convenience let us suppose that $r_1=q_1^{1/3}(1+\delta_1)$ and $r_2=1+\delta_2$, where $\delta_1,\delta_2\ll1$. Putting these values of $r_1$ and $r_2$ into the equations $r_1^{2}=(x+\mu)^{2}+y^2$ and $r_2^{2}=(x+\mu-1)^{2}+y^2$ and then solving with the rejection of second and higher order terms of $\delta_1$ and $\delta_2$, we get 

$x=\frac{q_1^{2/3}}{2}-\mu+(q_1^{2/3} \delta_1-\delta_2)$ and $y=\pm q_1^{1/3}\left(1-\frac{q_1^{2/3}}{4}+(2-q_1^{2/3}) \delta_1+\delta_2\right){}^{1/2}$.
We determine the values of $\delta_1$ and $\delta_2$ by putting the values of $x,y$ (after neglecting the terms containing $\delta_1,\delta_2\ll1$), $r_1=q_1^{1/3}(1+\delta_1)$ and $r_2=1+\delta_2$ into the equation (\ref{eq:wx}) and (\ref{eq:wy}) and solving them with the rejection of second and higher order terms of $\delta_1$ and $\delta_2$, we get \begin{eqnarray}
&&\delta_1=\frac{1}{3}\biggl[1-n^2 +2 c h \pi \frac{(b-a)}{a b}\frac{1}{(\mu^2 +q_1^{2/3}(1-\mu))^{}{3/2}}+\frac{3}{8} c h \pi \frac{\lgbya}{(\mu^2 +q_1^{2/3}(1-\mu))^2}\biggr], \label{eq:dl1} \\
&& \delta_2=\frac{1}{3(1+\frac{5}{2}A_2)}\biggl[1+\frac{3}{2}A_2-n^2 +2 c h \pi \frac{(b-a)}{a b}\frac{1}{(\mu^2 +q_1^{2/3}(1-\mu))^{}{3/2}}+\nonumber\\&&\frac{3}{8} c h \pi \frac{\lgbya}{(\mu^2 +q_1^{2/3}(1-\mu))^2}\biggr]. \label{eq:dl2} \end{eqnarray}
Hence, the coordinates of triangular points are :
\begin{eqnarray} \begin{cases}
x=\frac{q_1^{2/3}}{2}-\mu+(q_1^{2/3} \delta_1-\delta_2) \nonumber\\
y=\pm q_1^{1/3}\left(1-\frac{q_1^{2/3}}{4}+(2-q_1^{2/3}) \delta_1+\delta_2\right){}^{1/2} \end{cases} \label{eq:xy} \end{eqnarray} where $\delta_1$ and $\delta_2$ are given by equation (\ref{eq:dl1}) and (\ref{eq:dl2}) respectively.

Numerically the co-ordinates of triangular equilibrium points are $L_4=(0.366171, 0.641213)$, and $L_5=(0.366171, -0.641213)$ which are  obtained by using same parametric values as in collinear case into the equation (\ref{eq:fx}) and (\ref{eq:gy}) and solving them for $x$ and $y$. We also calculated $L_4=(0.447217, 0.515281)$ and $L_4=(0.447217, -0.515281)$ with the help of equation (\ref{eq:xy}) by evaluating $\delta_1$ and $\delta_2$ from the equations (\ref{eq:dl1}) and (\ref{eq:dl2}). The mean motion $(n)$ is calculated at disk's reference radius $(r)=0.99$ and found very similar results. We have seen that the  $L_{4(5)}$  are no longer remain triangular equilibrium points as they are in classical case.

\section{Zero Velocity Surface}
\label{sec:zvs}

Since, equation (\ref{eq:ji}) generally denotes the relation between coordinate and velocity of infinitesimal mass with respect to the rotating system, so by taking velocity term as zero in equation (\ref{eq:ji}) we have the equation of zero velocity surfaces as follows:
\begin{eqnarray}
&&2 \Omega=C \end{eqnarray} that is\begin{eqnarray}
\begin{cases}  
x^{2}+y^{2}-\frac{2(1-\mu)q_1}{r_1}-\frac{2\mu}{r_2}-\frac{\mu{A_2}}{r^{3}_2}+4 c h \pi \frac{(b-a)}{ab}\frac{1}{r}-\frac{7}{4}c h \pi \frac{\lgbya}{r^2}=C \nonumber\\
r_1=\sqrt{(x+\mu)^{2}+y^2 +z^2}\nonumber\\
r_2=\sqrt{(x+\mu-1)^{2}+y^2 +z^2}\nonumber\\
r=\sqrt{x^2 +y^2 +z^2}
\end{cases} \label{eq:zc} \end{eqnarray} 
Here, we try to know the approximate form of zero velocity surface[as in \cite{Moulton1914QB351.M92}] by analyzing the shape of the curves obtained by intersection of surface (\ref{eq:zc}) with the  $xy$-plane. The equation of curve is given by putting $z=0$ in (\ref{eq:zc}) which is as follows:
 \begin{figure}
 \plotone{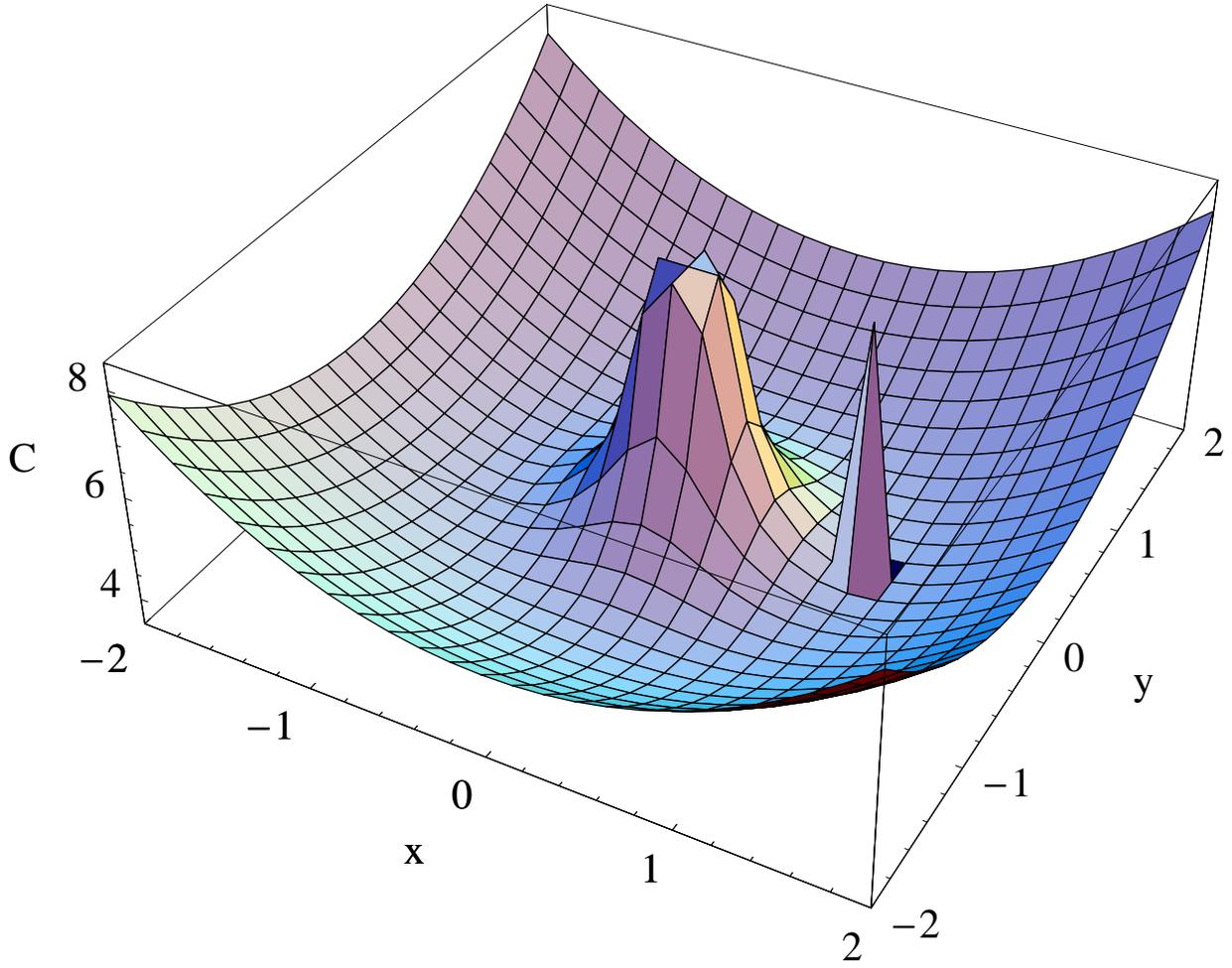}
 \caption{Zero velocity surface: $q_1=0.75$, $A_2=0.0025$, $a=1$, $b=1.5$ and  $r=0.99$.\label{zvc1_q75a25b15}}
 \end{figure}
\begin{eqnarray} 
&&x^{2}+y^{2}+\frac{2(1-\mu)q_1}{\sqrt{(x+\mu)^{2}+y^2}}+\frac{2\mu}{\sqrt{(x+\mu-1)^{2}+y^2}}+\frac{\mu{A_2}}{((x+\mu-1)^{2}+y^2){}^{3/2}}+\nonumber\\&&4 c h \pi\frac{(b-a)}{ab}\frac{1}{\sqrt{x^2 +y^2}}-\frac{7}{4}c h \pi\frac{\lgbya}{x^2 +y^2}=C 
\label{eq:zc1} \end{eqnarray} Here if we increases the value of $x$ and $y$, all terms, except first and second, on L.H.S. of (\ref{eq:zc1}), become negligible compare to first and second terms. Thus for the large value of $x$ and $y$, curves become very close to a circle of radius $\sqrt{C-\epsilon}$ having equation
\begin{eqnarray}
&&x^2 +y^2 =C-\frac{2(1-\mu)q_1}{\sqrt{(x+\mu)^{2}+y^2}}-\frac{2\mu}{\sqrt{(x+\mu-1)^{2}+y^2}}-\frac{\mu{A_2}}{((x+\mu-1)^{2}+y^2){}^{3/2}}-\nonumber\\&&4 c h \pi\frac{(b-a)}{ab}\frac{1}{\sqrt{x^2 +y^2}}+\frac{7}{4}c h \pi\frac{\lgbya}{x^2 +y^2}=C-\epsilon, \end{eqnarray} where, $\epsilon$ is very very small. For the small value of $x$ and $y$, the terms $x^2 +y^2$ will be negligible compare to remaining terms on L.H.S. of (\ref{eq:zc1}) and therefore, the curves become very close to an equipotential surface given bellow 
\begin{eqnarray} 
&&\frac{2(1-\mu)q_1}{\sqrt{(x+\mu)^{2}+y^2}}+\frac{2\mu}{\sqrt{(x+\mu-1)^{2}+y^2}}+\frac{\mu{A_2}}{((x+\mu-1)^{2}+y^2){}^{3/2}}+\nonumber\\&&4 c h \pi\frac{(b-a)}{ab}\frac{1}{\sqrt{x^2 +y^2}}-\frac{7}{4}c h \pi\frac{\lgbya}{x^2 +y^2}=C-x^{2}-y^{2}=C-\epsilon, \end{eqnarray} where $\epsilon$ is very very small.

The surface is given in figure  \ref{zvc1_q75a25b15} where two singular regions are shown by conic shapes.

\section{Stability of Equilibrium Points}
\label{sec:stb}
For the stability of equilibrium points $(x_e,y_e)$ let us assume a small change in its coordinate as $x=x_e+\xi$,\,$y=y_e+\eta$, where the displacements $\xi=P_1 e^{\lambda{t}}$,\,$\eta=P_2 e^{\lambda{t}}$ are very small, $P_1, P_2$ and $\lambda$ are parameters to be determined. Putting these coordinates into equations (\ref{eq:ux}) and (\ref{eq:vx}), we have the transformed equations of motion as follows[as in \cite{Murray2000ssd..book.....M}]
\begin{eqnarray}
\ddot{\xi}-2n\dot{\eta}&=&{\xi}{\Omega^0_{xx}}+{\eta}{\Omega^0_{xy}},  \\
\ddot{\eta}+2n\dot{\xi}&=&{\xi}{\Omega^0_{yx}}+{\eta}{\Omega^0_{yy}}, \end{eqnarray} where superfix $0$ denotes the corresponding value at equilibrium point. Now putting $\xi=P_1 e^{\lambda{t}}$,\,$\eta=P_2 e^{\lambda{t}}$ and simplifying we have
 \begin{eqnarray}
&&(\lambda^2-\Omega^0_{xx})P_1+(-2 n \lambda-\Omega^0_{xy})P_2=0, \\
&&(2 n \lambda-\Omega^0_{yx})P_1+(\lambda^2-\Omega^0_{yy})P_2=0, \end{eqnarray}
for nontrivial solution we have
\[
\begin{vmatrix}
\lambda^2-\Omega^0_{xx}& -2 n \lambda-\Omega^0_{xy} \\2 n \lambda-\Omega^0_{yx}& \lambda^2-\Omega^0_{yy}\\
\end{vmatrix}
=0.
\]
Simplifying above we have characteristic equation as follows:
\begin{eqnarray}
&&\lambda^4+(4 n^2-\Omega^0_{xx}-\Omega^0_{yy}){\lambda^2}+{(\Omega^0_{xx}}{\Omega^0_{yy}}-{\Omega^0}^2_{xy})=0.  \label{eq:ce} \end{eqnarray}
Now,\begin{eqnarray}
&&\Omega_{xx}=n^2-\frac{q_1(1-\mu)}{r^3_1}-\frac{\mu}{r^3_2}-\frac{3}{2}\frac{\mu {A_2}}{r^5_2}-2 c h \pi \frac{(b-a)}{a b}\frac{1}{r^3}\nonumber\\&&-\frac{3}{8} c h \pi {\lgbya}\frac{1}{r^4}+\frac{3q_1(1-\mu)(x+\mu)^2}{r^5_1}+\frac{3\mu(x+\mu-1)^2}{r^5_2}+\frac{15}{2}\frac{\mu {A_2}(x+\mu-1)^2}{r^7_2}\nonumber\\&&+6 c h \pi \frac{(b-a)}{a b}\frac{x^2}{r^5}+\frac{3}{2} c h \pi {\lgbya}\frac{x^2}{r^6},\\
&&\Omega_{yy}=n^2-\frac{q_1(1-\mu)}{r^3_1}-\frac{\mu}{r^3_2}-\frac{3}{2}\frac{\mu {A_2}}{r^5_2}-2 c h \pi \frac{(b-a)}{a b}\frac{1}{r^3}-\frac{3}{8} c h \pi {\lgbya}\frac{1}{r^4}+\nonumber\\&&\frac{3q_1(1-\mu) y^2}{r^5_1}+\frac{3\mu y^2}{r^5_2}+\frac{15}{2}\frac{\mu {A_2} y^2}{r^7_2}+6 c h \pi \frac{(b-a)}{a b}\frac{y^2}{r^5}+\frac{3}{2} c h \pi {\lgbya}\frac{2^2}{r^6},\\
&&\Omega_{xy}=\Omega_{yx}=\frac{3q_1(1-\mu)(x+\mu) y}{r^5_1}+\frac{3\mu(x+\mu-1) y}{r^5_2}+\frac{15}{2}\frac{\mu {A_2}(x+\mu-1) y}{r^7_2}\nonumber\\&&+6 c h \pi \frac{(b-a)}{a b}\frac{x y}{r^5}+\frac{3}{2} c h \pi {\lgbya}\frac{x y}{r^6}. \end{eqnarray}

In the following subsection  we examine the stability of equilibrium points.
\subsection{Stability of Collinear Equilibrium Points}
For the stability of collinear points $(y=0)$, we have 
\begin{eqnarray}
&&\Omega_{xx}=n^2+\frac{2q_1(1-\mu)}{r^3_1}+\frac{2\mu}{r^3_2}+\frac{6\mu {A_2}}{r^5_2}+4 c h \pi \frac{(b-a)}{a b}\frac{1}{r^3}\nonumber\\&&+\frac{9}{8} c h \pi {\lgbya}\frac{1}{r^4},\\ 
&&\Omega_{yy}=n^2-\frac{q_1(1-\mu)}{r^3_1}-\frac{\mu}{r^3_2}-\frac{3}{2}\frac{\mu {A_2}}{r^5_2}-2 c h \pi \frac{(b-a)}{a b}\frac{1}{r^3}\nonumber\\&&-\frac{3}{8} c h \pi {\lgbya}\frac{1}{r^4},\\ 
&&\Omega_{xy}=\Omega_{yx}=0, \end{eqnarray} also we have  $r_{1}=|x+\mu|$, $r_{2}=|x+\mu-1|$ and $r=|x|$.
To insure the stability of equilibrium points, we must have $\xi=P_1 e^{\lambda{t}}$,\,$\eta=P_2 e^{\lambda{t}}$ which can be expressed in a periodic functions. In other words the four roots $\lambda_1,\,\lambda_2,\,\lambda_3$ and $\lambda_4$ of the characteristic equation (\ref{eq:ce}) must be pure imaginary[as in \cite{Boccaletti1996thor.book.....B}] otherwise we have an unstable point. Hence for this we examine here nature of roots which is easily done by finding the sign of $\Omega^0_{xx}$ and  $\Omega^0_{yy}$ as follows:

In case(1)  when $1-\mu<x$, we have $\Omega^0_{xx}>0$ and $\Omega^0_{yy}<0$ and hence the discriminant $D=(4 n^2-\Omega^0_{xx}-\Omega^0_{yy})^2-4{\Omega^0_{xx}}{\Omega^0_{yy}}$ has positive sign. Therefore in this case, the two of the four roots $\lambda_1,\,\lambda_2,\,\lambda_3$ and $\lambda_4$ of the characteristic equation (\ref{eq:ce}) are complex conjugate pair and remaining two are pure imaginary conjugate pair. That is 
\begin{eqnarray}
&&\lambda_1=-\lambda_2=\sqrt{\frac{-\left(4 n^2-\Omega^0_{xx}-\Omega^0_{yy}\right)+\sqrt{\left(4 n^2-\Omega^0_{xx}-\Omega^0_{yy}\right)^2-4{\Omega^0_{xx}}{\Omega^0_{yy}}}}{2}}\nonumber\\&& \text{ is a complex conjugate pair and}\nonumber\\
&&\lambda_3=-\lambda_4=\sqrt{\frac{-\left(4 n^2-\Omega^0_{xx}-\Omega^0_{yy}\right)-\sqrt{\left(4 n^2-\Omega^0_{xx}-\Omega^0_{yy}\right)^2-4{\Omega^0_{xx}}{\Omega^0_{yy}}}}{2}}\nonumber\ \end{eqnarray}is a pure imaginary conjugate pair. Thus the condition of stability failed. Similarly, analyze the case (2)(i)  when $0\leq x<1-\mu$ we have  $\Omega^0_{xx}>0$ and $\Omega^0_{yy}>0$, in case (2)(ii) when, $-\mu<x<0$  we have $\Omega^0_{xx}>0$ and $\Omega^0_{yy}<0$ and in case (3) when, $x<-\mu$ we have $\Omega^0_{xx}>0$ and $\Omega^0_{yy}>0$. That is in each case we have a positive discriminant. Therefore, we conclude all the collinear equilibrium points are unstable. We examined the stability of collinear equilibrium points with respect to same parametric values setting above and conclude that generally $L_1$, $L_2$, $L_3$ and  $l_1$(new point) are unstable points but $L_2$ is stable for the disk's width ($a=1,\, 1.446\leq b\leq1.767$) also $L_3$ is stable for the disk's width ($a=1,\, 1\leq b\leq1.08$).

\subsection{Stability of Equilibrium Point $L_4$} 
\begin{figure}
 \plottwo{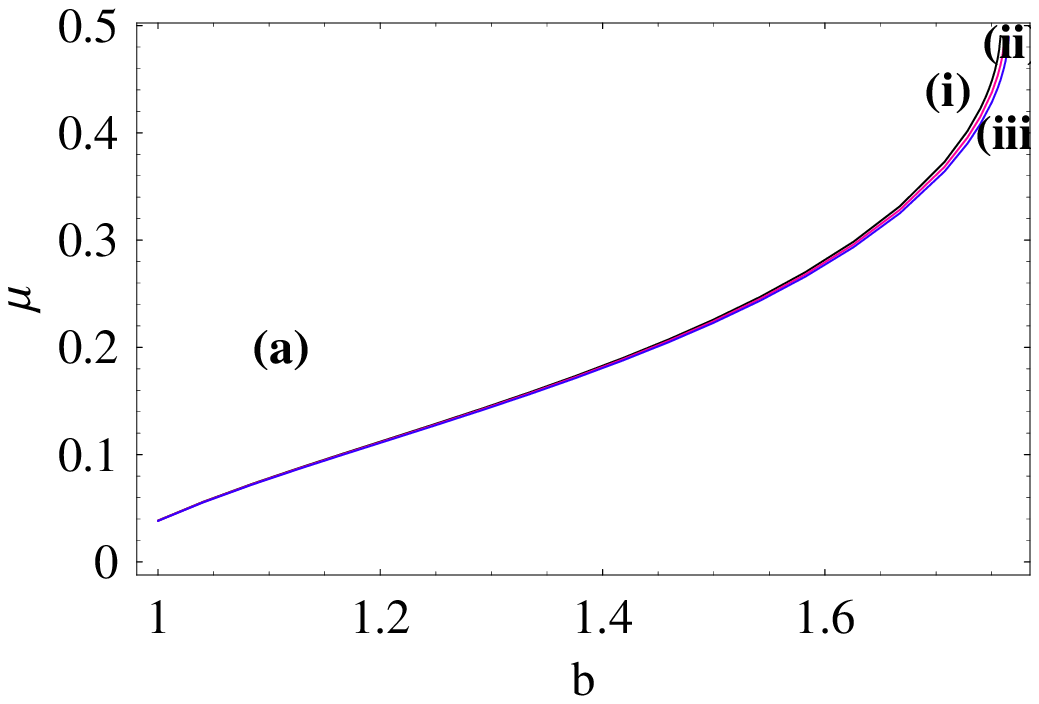}{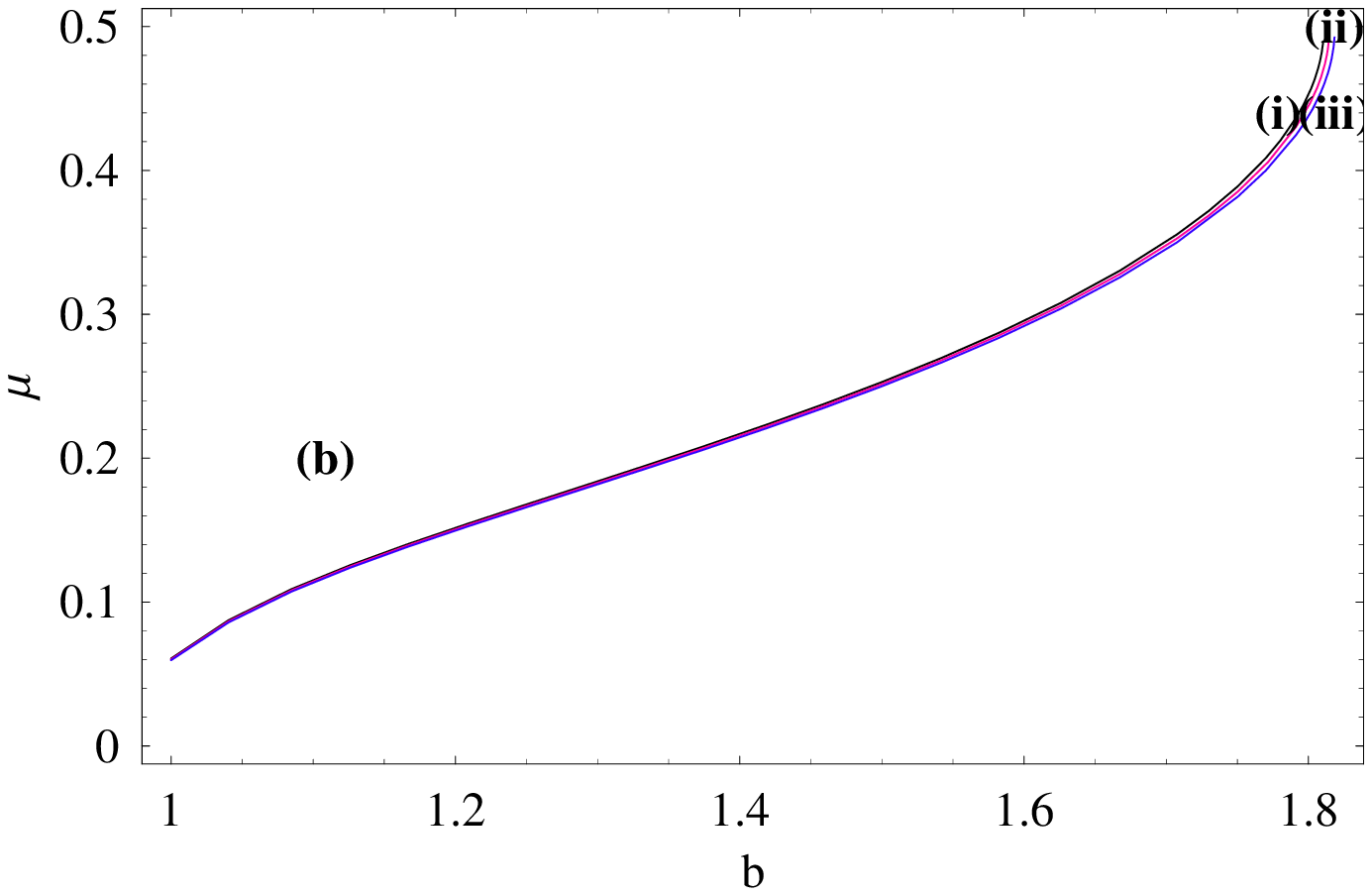}
 \caption{Variation of critical mass: (i) $A_2=0$, (ii) $A_2=0.0025$, (iii) $A_2=0.0050$. (a) at  $q_1=1, a=1, r=0.99$, (b) at  $q_1=0.75, a=1, r=0.99$.}\label{mAq1blable.eps}
\end{figure}
In case of $L_4$, we have \begin{eqnarray}
&&4 n^2-\Omega^0_{xx}-\Omega^0_{yy}=n^2-\frac{3\mu A_2}{r^5_{2,0}}-\frac{3}{8}c h \pi \lgbya\frac{1}{r^4_0},\\
&&{\Omega^0_{xx}}{\Omega^0_{yy}}-{\Omega^0}^2_{xy}=9 \mu(1-\mu) \gamma_0, \end{eqnarray} where \begin{eqnarray}
&&\gamma_0=y^2_0\biggl[\frac{q_1}{{r^5_{1,0}}{r^5_{2,0}}}\left(1+\frac{5 A_2}{2 r^2_{2,0}}\right)+\frac{q_1 \mu}{{r^5_{1,0}}{r^5_{0}}}\left(2\frac{b-a}{a b}+\frac{\log b-\log a}{2 r_0}\right)\nonumber\\&&+\frac{c h \pi(1-\mu)}{{r^5_{0}}{r^5_{2,0}}}\left(1+\frac{5 A_2}{2 r^2_{2,0}}\right)\left(2\frac{b-a}{a b}+\frac{\log b-\log a}{2 r_0}\right)\biggr] \end{eqnarray} The superfix $0$ indicate that values are at equilibrium points. Now, from the characteristic equation (~\ref{eq:ce}) we have\begin{eqnarray}
&&\lambda^2=\frac{-\left(4 n^2-\Omega^0_{xx}-\Omega^0_{yy}\right)\pm\sqrt{\left(4 n^2-\Omega^0_{xx}-\Omega^0_{yy}\right)^2-4\left({\Omega^0_{xx}}{\Omega^0_{yy}}-{\Omega^0}^2_{xy}\right)}}{2} \end{eqnarray} It is seen earlier that $\xi=P_1 e^{\lambda{t}}$,\,$\eta=P_2 e^{\lambda{t}}$ will be periodic and bounded if $\lambda_i$ ,$i=1,2,3,4$ is pure imaginary, i.e. $\lambda^2<0$. Therefore for stable solution we must have $\left(4 n^2-\Omega^0_{xx}-\Omega^0_{yy}\right)^2>4\left({\Omega^0_{xx}}{\Omega^0_{yy}}-{\Omega^0}^2_{xy}\right)$. In other words $\left(n^2-\frac{3\mu A_2}{r^5_{2,0}}-\frac{3}{8}c h \pi \lgbya\frac{1}{r^4_0}\right)^2>36\mu(1-\mu)\gamma_0$. For the classical values i.e.$(q_1=1, \,A_2=0, \,b=a, \,n=1, \,x=\frac{1}{2}-\mu, \,y=\pm\frac{\sqrt{3}}{2}, \,r_1=r_2=1)$, the above inequality reduces as $27\mu(1-\mu)<1$ which provides $\mu<0.0385209=\mu_c$ the value of critical mass in classical case. But in our case the value of critical mass for different values of $q_1$, $A_2$, and $b$ will be obtained by following graphs.

In figure \ref{mAq1blable.eps}, we plot $\mu$ Vs $b$, where (a) for $q_1=1$ and (b) for $q_1=0.75$ both containing three curves (i) $A_2=0$, (ii) $A_2=0.0025$, (iii) $A_2=0.0050$ have same nature i.e. as we increases the value of disk's width $b$, critical mass $\mu_c$ also increases. Initially these curves increases slowly but when $b>1.6$ they increase strictly as depicted in figure\ref{mAq1blable.eps}. These shows that when we increase the width of disk, stability region spans slowly then rapidly after $b=1.6$ upto $\mu_c=0.5$. 
\section{Conclusion}
\label{sec:conclude}
We have studied the dynamical properties of modified restricted three body problem and found that there exists a new equilibrium point in addition to five equilibrium points in classical problem. The zero velocity surface have obtained and examined the stability of equilibrium points. It is found that the $L_2$ and $L_3$ are stable for certain values of inner and outer radius of the disk and other collinear points are unstable while $L_4$ is conditionally stable upto the Routh's value of mass ratio. Thus stability region spans with the width of the disk. Hence, we conclude that the shape and size of the disk are  very significant for the motion of the bodies in space. 

\acknowledgements{This work is supported by the Department of Science and Technology, Govt. of India through the SERC-Fast Track Scheme for Young Scientist [SR/FTP/PS-121/2009]. First two  authors are also thankful to IUCAA Pune for partially financial support to visit library and to use computing facility.}

\end{document}